\documentclass[reqno,onecolumn,oneside]{paper}
\usepackage[english]{babel}
\usepackage{enumerate,amsmath,amsthm,amsfonts,pifont,slashed,amssymb,ifsym,mathrsfs,graphicx,graphics,yfonts,calligra,
latexsym,txfonts%,pxfonts
}
\usepackage[all]{xy}
\usepackage[sort,compress]{cite}

\newtheorem{theorem}{Theorem}[section]
\newtheorem{proposition}[theorem]{Proposition}
\newtheorem{lemma}[theorem]{Lemma}
\newtheorem{corollary}[theorem]{Corollary}
\newtheorem{conjecture}[theorem]{Conjecture}
\theoremstyle{definition}
\newtheorem{definition}[theorem]{Definition}
\newtheorem{exm}[theorem]{Example}

\theoremstyle{remark}
\newtheorem{remark}[theorem]{Remark}

\newcommand{\lcc}{\operatorname{lcc}}
\newcommand{\bydef}{\mathrel{\mathop:}=}

\newcommand{\op}{\operatorname{op}}

\newcommand{\id}{\operatorname{id}}

\newcommand{\Id}{\operatorname{Id}}
\newcommand{\Rad}{\operatorname{Rad}}

\newcommand{\supp}{\operatorname{supp}}

\newcommand{\BB}{\mathcal{B}\!\operatorname{oole}}
\renewcommand{\O}{\Omega}
\newcommand{\Set}{\mathcal{S}\!\operatorname{et}}
\newcommand{\fwp}{\mathscr{F}}

\newcommand{\Q}{\mathbb{Q}}
\renewcommand{\P}{\mathbb{P}}

\newcommand{\MV}{\mathcal{MV}}

\newcommand{\cat}{\mathcal}

\newcommand{\B}{\operatorname{B}}
\newcommand{\sk}{\operatorname{Sk}}
\newcommand{\clop}{\operatorname{Clop}}
\newcommand{\Max}{\operatorname{Max}}

\newcommand{\R}{\mathbb{R}}
\newcommand{\N}{\mathbb{N}}
\newcommand{\C}{\mathbb{C}}
\newcommand{\Z}{\mathbb{Z}}

\renewcommand{\wp}{\mathscr{P}}

\newcommand{\TMV}{{}^{\operatorname{MV}}\!\mathcal{T}\!\!\operatorname{op}}
\renewcommand{\Top}{\mathcal{T}\!\!\operatorname{op}}
\newcommand{\SMV}{{}^{\operatorname{MV}}\!\!\mathcal{S}\!\operatorname{tone}}
\newcommand{\MVs}{\mathcal{MV}^{\operatorname{ss}}}
\newcommand{\MVsfc}{\mathcal{MV}^{\operatorname{sfc}}}
\newcommand{\MVlcc}{\mathcal{MV}^{\operatorname{lcc}}}
\newcommand{\ST}{\mathcal{S}\!\operatorname{tone}}
\newcommand{\app}{\approx^{\operatorname{op}}}

\renewcommand{\phi}{\varphi}

\newcommand{\G}{\mathbf G}
\newcommand{\cou}{{^\leftarrow}}
\newcommand{\fcou}{{^{\rotatebox[origin=c]{180}{$\rightsquigarrow$}}}}
\newcommand{\restr}{\upharpoonright}
\newcommand{\la}{\langle}
\newcommand{\ra}{\rangle}

\newcommand{\lto}{\longrightarrow}
\newcommand{\To}{\Longrightarrow}
\newcommand{\lmapsto}{\longmapsto}

\newcommand{\ul}{\underline}

\newcommand{\0}{\mathbf{0}}
\renewcommand{\1}{\mathbf{1}}
\renewcommand{\2}{\mathbf{2}}
\def\amslatex\slash{{\protect\AmS-\protect\LaTeX}}

\begin{document} 

\title{An extension of Stone Duality to fuzzy topologies and MV-algebras}

\author{Ciro Russo\thanks{This work was carried out within the IRSES project MaToMUVI, funded by the EU 7th Framework Programme.}} 

\institution{Departamento de Matem\'atica \\ Instituto de de Matem\'atica \\ Universidade Federal da Bahia, Brazil \\ \small{\texttt{ciro.russo@ufba.br}}}

\maketitle
\today

\begin{abstract}
In this paper we introduce the concept of \emph{MV-topology}, a special class of fuzzy topological spaces, and prove a proper extension of Stone Duality to the categories of \emph{limit cut complete MV-algebras} and \emph{Stone MV-spaces}, namely, zero-dimensional compact Hausdorff MV-topological spaces. Then we describe the object class of limit cut complete MV-algebras, and show that any semisimple MV-algebra has a limit cut completion, namely, a minimum limit cut complete extension. Last, we compose our duality with other known ones, thus obtaining new categorical equivalences and dualities involving categories of MV-algebras.
\end{abstract}

\section{Introduction}

The concept of \emph{fuzzy topology} was introduced a few years after Zadeh's famous paper on fuzzy sets \cite{zad}, and its study has been pursued for many years (see, for instance, \cite{chal,hoh1,hoh2,hoh3,low,maru,rod1,rod2,sto1,sto2}). In defining a fuzzy topological space on a set $X$ a fundamental role is played by the structure used to represent the ``fuzzy powerset'' of $X$, i.e., the fuzzy version of the Boolean algebra $\2^X$. According to the original definition of fuzzy set, one may find natural to consider $[0,1]^X$ as the fuzzy powerset of $X$. As a matter of fact, most of the authors in this area approached fuzzy topology using either arbitrary lattice-valued fuzzy subsets or $[0,1]^X$ with its natural lattice structure. However, fuzzy topological spaces using $[0,1]^X$ equipped with a richer algebraic structure (e. g., continuous or left-continuous t-norms \cite{haj}) have been considered in the literature. In our opinion, looking at the crisp and fuzzy powersets of a given set $X$ as, respectively, $\2^X$ and $[0,1]^X$, it is undoubtable that the structure of MV-algebra \cite{cha1} of the latter is the one that best succeeds in preserving many properties of symmetry that are inborn qualities of Boolean algebras.

On the other hand, the duality theory for MV-algebras boasts a rather wide interest among researchers in the area \cite{cidumu,latpet,marspa,gvgm, chk,cima,dnl,cimu}, including some of the most prominent ones, but --- quite surprisingly, indeed --- the only relevant work connecting MV-algebras and fuzzy topologies via a duality is, to the best of our knowledge, a paper by Maruyama \cite{maru} published in 2010. Such a circumstance is even more curious if we consider that a Stone-type representation theorem for semisimple MV-algebras was published in 1986 \cite{bel} but probably foreseen since right after the pioneering work of Chang \cite{cha1}.

In this paper we propose a concept of fuzzy topological space which is a natural generalization of classical topology with the use of MV-algebras. Moreover, even if the present paper is more of an algebraic and categorical nature, a strong motivation for the introduction of such fuzzy topologies comes from the area of Mathematical Morphology \cite{hei,mat} and its connection to quantales and idempotent semirings discussed in \cite{rus}. Indeed, as binary digital images are examples of crisp subsets of a given set, greyscale images are a prototypical example of fuzzy subsets. So, since mathematical morphological operators, for binary images, are designed with the aim of individualizing relevant topological properties of the images, it is clear that classical topology can hardly give satisfactory information in the case of greyscale images.

Our aim is to use MV-algebras as a framework for fuzzy topology which, on the one hand, is sufficiently rich and complex and, on the other hand, reflects (up to a suitable reformulation) as many properties of classical topology as possible. For this reason we introduce the concept of \emph{MV-topology}, a generalization of general topology whose main features can be summarized as follows.
\begin{itemize}
\item The Boolean algebra of subsets of the universe is replaced by the MV-algebra of ($[0,1]$-valued) fuzzy subsets.
\item Classical topological spaces are examples of MV-topological spaces.
\item The algebraic structure of the family of open (fuzzy) subsets has a quantale reduct $\la \O, \bigvee, \oplus\ra$, which replaces the classical sup-lattice $\la \O, \bigvee\ra$, and an idempotent semiring one $\la \O, \wedge, \odot, \1\ra$ in place of the meet-semilattice $\la\O, \wedge, \1\ra$. Moreover, the lattice reduct $\la \O, \bigvee, \wedge\ra$ maintains the property of being a frame.
\item The MV-algebraic negation $^*$ is, in the aforementioned classes of algebras, an isomorphism between the various structures of open subsets and the corresponding ones of closed subsets.
\item A classical topology is canonically associated to each MV-topology. It is called the \emph{skeleton topology} and is obtained simply by restricting the family of open subsets to the crisp ones.
\end{itemize}

The main results of the paper are proved in Sections \ref{duality}, \ref{lccsec}, and \ref{sec:comb}.

In particular, we show (Section \ref{duality}) an extension of Stone Duality between Boolean algebras and Stone spaces to, respectively, the category of \emph{limit cut complete MV-algebras}, namely, the full subcategory of $\MV$ whose objects are algebras which contain the suprema of certain cuts, and a suitable category of MV-topologies, whose objects are the natural MV-version of Stone (or Boolean) spaces --- called \emph{Stone MV-spaces}. Such an extension is ``proper'' in the sense that its restriction to, respectively, Boolean algebras and Stone spaces --- which are full subcategories of the ones involved in the duality --- yields the classical well-known duality, up to a trivial reformulation in terms of maximal ideals instead of ultrafilters. In Section \ref{lccsec} we shall describe limit cut complete MV-algebras and characterize the clopen algebras of strongly compact Stone MV-spaces. Last, in Section \ref{sec:comb} we shall connect our duality theorem with other known dualities for classes of MV-algebras, thus obtaining new categorical equivalences and dualities.

We refer the reader to the reference books \cite{mvbook,mubook} for all the necessary notions and results on MV-algebras not explicitly reported here.

\section{MV-topologies}
\label{sec:top}

Throughout the paper, both crisp and fuzzy subsets of a given set will be identified with their membership functions and usually denoted by lower case latin or greek letters. In particular, for any set $X$, we shall use also $\1$ and $\0$ for denoting, respectively, $X$ and $\varnothing$. In some cases, we shall use capital letters in order to emphasize that the subset we are dealing with is crisp.

We remark that an MV-topological space is basically a special fuzzy topological space in the sense of C. L. Chang \cite{chal}. Moreover, most of the definitions and results of the present section and of Section \ref{sec:base} are simple adaptations of the corresponding ones of the aforementioned work to the present context or directly derivable from the same work or from the results presented in the papers \cite{hoh1,hoh2,hoh3,low,maru,rod1,rod2,sto1,sto2} that we already cited in the introduction. 

\begin{definition}\label{mvtop}
Let $X$ be a set, $A$ the MV-algebra $[0,1]^X$ and $\O \subseteq A$. We say that $\la X, \O\ra$ is an \emph{MV-topological space} if $\O$ is a subuniverse both of the quantale $\la [0,1]^X, \bigvee, \oplus\ra$ and of the semiring $\la [0,1]^X, \wedge, \odot, \1\ra$. More explicitly, $\la X, \O\ra$ is an MV-topological space if
\begin{enumerate}[(i)]
\item $\0, \1 \in \O$,
\item for any family $\{o_i\}_{i \in I}$ of elements of $\O$, $\bigvee_{i \in I} o_i \in \O$,
\end{enumerate}
and, for all $o_1, o_2 \in \O$,
\begin{enumerate}[(i)]
\setcounter{enumi}{2}
\item $o_1 \odot o_2 \in \O$,
\item $o_1 \oplus o_2 \in \O$,
\item $o_1 \wedge o_2 \in \O$.
\end{enumerate}
$\O$ is also called an \emph{MV-topology} on $X$ and the elements of $\O$ are the \emph{open MV-subsets} of $X$. The set $\Xi = \{o^* \mid o \in \O\}$ is easily seen to be a subquantale of $\la [0,1]^X, \bigwedge, \odot\ra$ (where $\bigwedge$ has to be considered as the join w.r.t. to the dual order $\geq$ on  $[0,1]^X$) and a subsemiring of $\la [0,1]^X, \vee, \oplus, \0\ra$, i.e., it verifies the following properties:
\begin{itemize}
\item[$-$] $\0, \1 \in \Xi$,
\item[$-$] for any family $\{c_i\}_{i \in I}$ of elements of $\Xi$, $\bigwedge_{i \in I} c_i \in \Xi$,
\item[$-$] for all $c_1, c_2 \in \Xi$, $c_1 \odot c_2, c_1 \oplus c_2, c_1 \vee c_2 \in \Xi$.
\end{itemize}
The elements of $\Xi$ are called the \emph{closed MV-subsets} of $X$.
\end{definition}

\begin{proposition}\label{sub}
Let $\la X, \O \ra$ be an MV-topological space. For any subset $Y$ of $X$, the pair $\la Y, \O_Y \ra$, where $\O_Y \bydef \{o_{\restr Y} \mid o \in \O\}$, is an MV-topology on $Y$.
\end{proposition}
\begin{proof}
Trivial.
\end{proof}

\begin{definition}\label{subspace}
For any subset $Y$ of $X$, the pair $\la Y, \O_Y\ra$ is called an \emph{MV-subspace} of $\la X, \O\ra$.
\end{definition}

\begin{exm}\label{ex1}
\begin{enumerate}[(a)]
\item $\la X, \{\0,\1\}\ra$ and $\la X, [0,1]^X \ra$ are MV-topological spaces.
\item Any topology is an MV-topology.
\item Let $d: X \lto [0,+\infty[$ be a distance function on $X$. For any fuzzy point $\alpha$ of $X$, with support $x$, and any positive real number $r$, we define the \emph{open ball} of center $\alpha$ and radius $r$ as the fuzzy set $\beta_r(\alpha)$ identified by the membership function $\beta_r(\alpha)(y) = \left\{\begin{array}{ll} \alpha(x) & \textrm{if } d(x,y) < r \\ 0 & \textrm{if } d(x,y) \geq r \end{array}\right.$. Analogously, the \emph{closed ball} $\beta_r[\alpha]$ of center $\alpha$ and radius $r$ has membership function $\beta_r[\alpha](y) = \left\{\begin{array}{ll} \alpha(x) & \textrm{if } d(x,y) \leq r \\ 0 & \textrm{if } d(x,y) > r \end{array}\right.$. It is immediate to verify that the fuzzy subsets of $X$ that are join of a family of open balls is an MV-topology on $X$ that is said to be \emph{induced} by $d$. This example can be found also in \cite{low}.
\end{enumerate}
\end{exm}

\begin{definition}\label{skeleton}
If $\la X, \O\ra$ is an MV-topology, then $\la X, \B(\O)\ra$ --- where $\B(\O) \bydef \O \cap \{0,1\}^X = \O \cap \B([0,1]^X)$ --- is both an MV-topology and a topology in the classical sense. The topological space $\la X, \B(\O)\ra$ will be called the \emph{skeleton space} of $\la X, \O\ra$.
\end{definition}

Observe that the skeleton space of a given MV-topological one can be equivalently defined by
$$\B(\O) = \{\Delta \circ \alpha \mid \alpha \in \O\},$$
where $\Delta$ is the so-called \emph{Baaz delta} operator \cite{baaz}, i.e.,
$$\Delta: x \in [0,1] \mapsto \left\{\begin{array}{l} 1 \text{ if x = 1} \\ 0 \text{ if x < 1} \end{array}\right. \in \{0,1\}.$$
$\Delta$, besides being a monotonic map, is a monoid homomorphism between $\la [0,1], \odot, 1 \ra$ and $\la \{0,1\}, \wedge, 1 \ra$. Therefore the equivalence of the two definitions follows from the fact that MV-topologies are closed under $\odot$ while classical ones are closed under $\wedge$.

Let $X$ and $Y$ be sets. Any function $f: X \lto Y$ naturally defines a map
\begin{equation}\label{muf}
\begin{array}{cccc}
f\fcou: & [0,1]^Y & \lto 		& [0,1]^X \\
		 & \alpha  & \lmapsto & \alpha \circ f.
\end{array}
\end{equation}
Obviously $f\fcou(\0) = \0$; moreover, if $\alpha, \beta \in [0,1]^Y$, for all $x \in X$ we have $f\fcou(\alpha \oplus \beta)(x) = (\alpha \oplus \beta)(f(x)) = \alpha(f(x)) \oplus \beta(f(x)) = f\fcou(\alpha)(x) \oplus f\fcou(\beta)(x)$ and, analogously, $f\fcou(\alpha^*) = f\fcou(\alpha)^*$. Then $f\fcou$ is an MV-algebra homomorphism and we shall call it the \emph{MV-preimage} of $f$. The reason of such a name is essentially the fact that $f\fcou$ can be seen as the preimage, via $f$, of the fuzzy subsets of $Y$. From a categorical viewpoint, once denoted by $\Set$, $\BB$ and $\MV$ the categories of sets, Boolean algebras, and MV-algebras respectively (with the obvious morphisms), there exist two contravariant functors $\wp: \Set \lto \BB^{\op}$ and $\fwp: \Set \lto \MV^{\op}$ sending each map $f: X \lto Y$, respectively, to the Boolean algebra homomorphism $f\cou: \wp(Y) \lto \wp(X)$ and to the MV-homomorphism $f\fcou: [0,1]^Y \lto [0,1]^X$.

Moreover, for any map $f: X \lto Y$ we define also a map $f^\to: [0,1]^X \lto [0,1]^Y$ by setting, for all $\alpha \in [0,1]^X$ and for all $y \in Y$,
\begin{equation}\label{ovf}
f^\to(\alpha)(y) = \bigvee_{f(x) = y} \alpha(x).
\end{equation}
Clearly, if $y \notin f[X]$, $f^\to(\alpha)(y) = \bigvee \varnothing = \0$ for any $\alpha \in [0,1]^X$.

\begin{definition}\cite{chal}\label{cont}
Let $\la X, \O_X \ra$ and $\la Y, \O_Y \ra$ be two MV-topological spaces. A map $f: X \lto Y$ is said to be
\begin{itemize}
\item \emph{continuous} if $f\fcou[\O_Y] \subseteq \O_X$,
\item \emph{open} if $f^\to(o) \in \O_Y$ for all $o \in \O_X$,
\item \emph{closed} if $f^\to(c) \in \Xi_Y$ for all $c \in \Xi_X$ 
\item an \emph{MV-homeomorphism} if it is bijective and both $f$ and $f^{-1}$ are continuous.
\end{itemize}
\end{definition}
We can use the same words of the classical case because, as it is trivial to verify, if a map between two classical topological spaces is continuous, open, or closed in the sense of the definition above, then it has the same property in the classical sense. 

Continuity, as in Definition \ref{cont}, is equivalent to $f\fcou[\Xi_Y] \subseteq \Xi_X$. Indeed, since $f\fcou: [0,1]^Y \lto [0,1]^X$ is an MV-algebra homomorphism, it preserves $^*$; therefore, for any closed set $c$ of $Y$, $c^*$ is an open set, hence $f\fcou(c^*) = f\fcou(c)^* \in \O_X$ implies $f\fcou(c) \in \Xi_X$. In a completely analogous way, it can be proved that $f\fcou[\Xi_Y] \subseteq \Xi_X$ implies continuity in the sense of the previous definition.

Moreover, it is absolutely obvious that, if $\la X, \O_X \ra$ and $\la Y, \O_Y \ra$ are two MV-spaces, and $f: X \lto Y$ is a continuous function between them, then $f$ is also a continuous map between the two skeleton spaces $\la X, \B(\O_X)\ra$ and $\la Y, \B(\O_Y)\ra$.

\section{Bases, compactness and separation axioms}
\label{sec:base}

In the present section we give the necessary definitions and show some preliminary results in order to prove the extension of Stone Duality. As the reader will notice, the concepts we are going to introduce are direct and natural (and quite obvious, indeed) generalizations of the corresponding ones in classical topology. Actually, some of the following definitions and results are either already present, or plainly adapted from similar ones, in the theory of Fuzzy Topology. In those cases, we shall give a suitable bibliographical reference.

In order to build a comprehensive theory of MV-topologies, many further material needs to be defined and investigated; nonetheless, as we already mentioned, here we focus our attention to Stone Duality. Therefore we do not intend to introduce notions that shall not be of any utility in this particular paper, leaving such further insights for future works.

\begin{definition}\cite{war}\label{base}
As in classical topology, we say that, given an MV-topological space $\tau = \la X, \O\ra$, a subset $\Theta$ of $[0,1]^X$ is called a \emph{base} for $\tau$ if $\Theta \subseteq \O$ and every open set of $\tau$ is a join of elements of $\Theta$. 
\end{definition}

\begin{lemma}\label{bascont}
Let $\tau = \la X, \O_X\ra$ and $\tau' = \la Y, \O_Y\ra$ be two MV-topological spaces and let $\Theta$ be a base for $\tau'$. A map  $f: X \lto Y$ is continuous if and only if $f\fcou[\Theta] \subseteq \O_X$.
\end{lemma}
\begin{proof}
One implication is trivial, since $\Theta$ is a family of open sets. Conversely, assuming that $f\fcou[\Theta] \subseteq \O_X$, let $o = \bigvee \Gamma$, with $\Gamma \subseteq \Theta$, be any open set of $\tau'$. As we observed, $f\fcou$ is an MV-algebra homomorphism, hence $f\fcou(o) = f\fcou\left(\bigvee \Gamma\right) = \bigvee f\fcou[\Gamma]$, i.e. $f\fcou(o)$ is the join of open sets of $\tau$ and, therefore, open itself.
\end{proof}

A \emph{covering} of $X$ is any subset $\Gamma$ of $[0,1]^X$ such that $\bigvee \Gamma = \1$ \cite{chal}, while an \emph{additive covering} ($\oplus$-covering, for short) is a finite family $\{\alpha_i\}_{i=1}^n$ of elements of $[0,1]^X$, $n < \omega$, such that $\alpha_1 \oplus \cdots \oplus \alpha_n = \1$. It is worthwhile remarking that we used the expression ``finite family'' in order to include the possibility for such a family to have repetitions. In other words, an additive covering is a finite subset $\{\alpha_1, \ldots, \alpha_k\}$ of $[0,1]^X$, along with natural numbers $n_1, \ldots, n_k$, such that $n_1\alpha_1 \oplus \cdots \oplus n_k \alpha_k = \1$. 

\begin{proposition}\label{cov}
For any set $X$, any covering of fuzzy subsets of $X$ which is closed under $\oplus$, $\odot$, and $\wedge$ is a base for an MV-topology on $X$. 
\end{proposition}
\begin{proof}
Let $\Gamma \subseteq [0,1]^X$ be a covering closed under $\oplus$, $\odot$, and $\wedge$, and let $\O = \{\bigvee G \mid G \subseteq \Gamma\}$. We have $\1 \in \O$, by definition of covering, and $\0 = \bigvee \varnothing \in \O$.

On the other hand, $\O$ is trivially closed under arbitray joins and $\odot$, $\oplus$, and $\wedge$ distribute over any existing join. Then, given $o_1, o_2 \in \O$, $o_1 = \bigvee_{i \in I} \alpha_i$ and $o_2 = \bigvee_{j \in J} \beta_j$, with $\{\alpha_i\}_{i \in I}, \{\beta_j\}_{j \in J} \subseteq \Gamma$, whence
$$o_1 \bullet o_2 = \left(\bigvee_{i \in I} \alpha_i\right) \bullet \left(\bigvee_{j \in J} \beta_j\right) = \bigvee_{i \in I} \left(\alpha_i \bullet \bigvee_{j \in J} \beta_j\right) = \bigvee_{i \in I} \bigvee_{j \in J} (\alpha_i \bullet \beta_j),$$
for $\bullet \in \{\oplus, \odot, \wedge\}$. So $\O$ verifies Definition \ref{mvtop}, i.e. it is an MV-topology, and $\Gamma$ is a base for it.
\end{proof}

%If $\la X, \O \ra$ is an MV-topological space, $\O$ naturally defines \emph{MV-interior} and \emph{MV-closure} operators:
%\begin{itemize}
%\item $Y^\circ \bydef \bigvee \{o \in \O \mid o \leq Y\}$,
%\item $\ov Y := \bigwedge \{c \in \Xi \mid Y \leq c\}$,
%\end{itemize}
%for all $Y \in [0,1]^X$. Obviously $Y^\circ \in \O$ and $\ov Y \in \Xi$ for all $Y \in [0,1]^X$. Moreover,
%%%%as in the classical case such operators are frame homomorphisms from $\la \wp(X), \bigcup, \cap, \varnothing, X\ra$ to $\O$ and from $\la \wp(X), \bigcap, \cup, X, \varnothing\ra$ to $\Xi$ respectively, for MV-topologies they are quantale homomorphisms from $\la [0,1]^X, \bigvee, \odot, \0, \1\ra$ to $\la \O, \bigvee, \odot, \0, \1\ra$ and from $\la [0,1]^X, \bigwedge, \oplus, \1, \0\ra$ to $\la \Xi, \bigwedge, \oplus, \1, \0\ra$ respectively. It is worth noticing that the MV-interior operator does not preserve $\oplus$ and, dually, the MV-closure does not preserve $\odot$. What 
%we have, for $\bullet \in \{\oplus, \odot\}$, $Y_1^\circ \bullet Y_2^\circ \leq (Y_1 \bullet Y_2)^\circ$ and $\ov{(Y_1 \bullet Y_2)} \leq \ov Y_1 \bullet \ov Y_2$.

The presence of strong and weak conjunctions and disjunction, in the structure of open sets of an MV-topology, naturally suggests different fuzzy versions (weaker or stronger) of most of the classical topological concepts (separation axioms, compactness etc.). However, we shall limit our attention to the ones that serve the scope of this paper, namely \emph{compactness} and \emph{Hausdorff (or $T_2$) separation axiom}. 

\begin{definition}\label{compact}
An MV-topological space $\la X, \O \ra$ is said to be \emph{compact} if any open covering of $X$ contains an additive covering; it is called \emph{strongly compact} if any open covering contains a finite covering.\footnote{What we call strong compactness here is called simply compactness in the theory of lattice-valued fuzzy topologies \cite{chal}.}
\end{definition}

It is obvious that strong compactness implies compactness and, since the operations $\oplus$ and $\vee$ coincide on Boolean elements of MV-algebras, in the case of topologies of crisp subsets the two notions collapse to the classical one. For the same reason, it is evident as well that the skeleton spaces of both compact and strongly compact MV-spaces are compact. The following example shows that compactness does not imply strong compactness, i.e., they are not equivalent.

\begin{exm}
Let $X$ be a non-empty set and $\O$ the set of all constant fuzzy subsets of $X$, which is clearly an MV-topology. For each $r \in [0,1]$, let $o_r$ be the fuzzy set constantly equal to $r$. Then, for any family $\{r_i\}_{i \in I} \subseteq [0,1)$ such that $\bigvee_{i \in I} r_i = 1$, the set $\{o_{r_i} \mid i \in I\}$ is an open covering and all the coverings not containing $\1$ are of this form. On the other hand, all of such coverings do not contain finite coverings but do include additive ones. 
\end{exm}

\begin{lemma}\label{closcomp}
A closed subspace $\la Y, \O_Y\ra$ of a compact (respectively: strongly compact) space $\la X, \O\ra$ is compact (resp.: strongly compact).
\end{lemma}
\begin{proof}
Since $Y$ is a subspace, in particular it is a crisp subset of $X$ and, therefore, all of its open sets are of the form $o_{\restr Y}$ with $o \in \O$. So let $\{o_i\}_{i \in I} \subseteq \O$ such that $\bigvee_{i \in I} o_i \geq Y$. Since $Y$ is closed, $Y^*$ is open and $\{o_i\}_{i \in I} \cup \{Y^*\}$ is an open covering of $X$. By compactness of $X$, there exists a finite family $\{o_j\}_{j=1}^n$ of elements of $\{o_i\}_{i \in I}$ such that $o_1 \oplus \cdots \oplus o_n \oplus Y^* = X$. Then, since $Y \wedge Y^* = \0$, we have (with a slight abuse of notation) $Y = Y \wedge (o_1 \oplus \cdots \oplus o_n) = (Y \wedge o_1) \oplus \cdots \oplus (Y \wedge o_n)$, the latter equality easily following from the properties of Boolean elements of MV-algebras, whence $Y$ is compact.

The case of strong compactness is completely analogous.
\end{proof}

\begin{definition}\label{t2}
Let $\tau = \la X, \O \ra$ be an MV-topological space. $X$ is called a \emph{Hausdorff} (or \emph{separated}) \emph{space} if, for all $x \neq y \in X$, there exist $o_x, o_y \in \O$ such that
\begin{enumerate}[(i)]
\item $o_x(x) = o_y(y) = 1$,
\item $o_x \wedge o_y = \0$.
\end{enumerate}
\end{definition}

\begin{remark}\label{t2r}
It is important to observe here that there is no interesting ``weak'' version of the above definition, since it is immediate to verify that Definition \ref{t2} is equivalent to the following:

for all $x \neq y \in X$, there exist $o_x', o_y' \in \O$ verifying
\begin{enumerate}
\item[(i)] $o_x'(x) = o_y'(y) = 1$,
\item[(ii')] $o_x \odot o_y = \0$.
\end{enumerate}
Indeed, overlooking the trivial implication, assume there such two open sets $o_x'$ and $o_y'$ exist, and set $o_x = o_x'^2$ and $o_y = o_y'^2$. Then, by the quasi-equation $x \odot y = 0 \To x^2 \wedge y^2 = 0$ (which holds in every MV-algebra), $o_x$ and $o_y$ satisfy Definition \ref{t2}.
\end{remark}

As for compactness, Definition \ref{t2} coincide with the classical $T_2$ property on crisp topologies and implies that the corresponding skeleton space is Hausdorff in the classical sense.  

The following result is obvious.
\begin{lemma}\label{xclos}
If $\la X, \O\ra$ is an Hausdorff space, then all crisp singletons of $X$ are closed. 
\end{lemma}

\section{The extension of Stone Duality}
\label{duality}

In this section we shall prove that Stone Duality can be extended to a class of semisimple MV-algebras and compact separated MV-topologies having a base of clopens. Before proving the duality theorem, we recall the definition of simple and semisimple MV-algebra along with a well-known representation theorem for the latter.

\begin{definition}{ss}
An MV-algebra $A$ is called \emph{simple} if its only proper ideal is $\{0\}$. $A$ is called \emph{semisimple} if it is a subdirect product of simple MV-algebras.
\end{definition}

It is well-known (see, for instance, \cite{mvbook}) that an MV-algebra $A$ is simple if and only if it is isomorphic to a subalgebra of $[0,1]$, and that $A$ is semisimple if and only if the \emph{radical of $A$}, $\Rad A$, i.e. the intersection of all maximal ideals of $A$, is $\{0\}$.

\begin{theorem}\cite{bel,cha1,cha2}\label{belrepr}
For any set $X$, the MV-algebra $[0,1]^X$ and all of its subalgebras are semisimple. Moreover, up to isomorphisms, all the semisimple MV-algebras are of this type. More precisely, every semisimple MV-algebra can be embedded in the MV-algebra of fuzzy subsets $[0,1]^{\Max A}$ of the maximal spectrum of $A$.
\end{theorem}
The proof of the first part of Theorem \ref{belrepr} is rather obvious. Before proving our main theorem, it is useful to briefly sketch the proof of the fact that any semisimple MV-algebra is embeddable in $[0,1]^{\Max A}$.
\begin{proof}
\textit{(Sketch)} For any maximal ideal $M$ the quotient algebra $A/M$ is a simple MV-algebra and, therefore, an Archimedean MV-chain. Then $A/M$ is isomorphic to a subalgebra of $[0,1]$ and we have this situation:
\begin{itemize}
\item for each $M \in \Max A$, there is the natural projection $\pi_M: A \lto A/M$;
\item for each $M \in \Max A$, there exists a unique embedding $\iota_M: A/M \lto [0,1]$;
\item the embedding $\iota: A \lto [0,1]^{\Max A}$ associates, to each $a \in A$, the fuzzy subset $\widehat a$ of $\Max A$ defined by $\widehat a(M) = \iota_M(\pi_M(a)) = \iota_M(a/M)$ for all $M \in \Max A$. 
\end{itemize}
\end{proof}
It is important to notice that the above construction is possible for any MV-algebra $A$ with the only difference (important, indeed) that the homomorphism $\iota$ is not injective if $A$ is not semisimple for the simple reason that $\ker \iota$ always coincides with $\Rad A$.

We will now recall some well-known properties of ideals of MV-algebras which shall be used in the subsequent proofs.
\begin{proposition}\cite{mvbook}\label{maxch}\label{properid}\label{classch}
Let $A$ be an MV-algebra, $I \in \Id(A)$, and $S \subseteq A$. Then the following hold.
\begin{enumerate}
\item[(i)] $I$ is maximal if and only if, for any $a \in A$, $a \notin I$ implies that there exists $n < \omega$ such that $(a^*)^n \in I$.
\item[(ii)] For all $a \in A$, $a/I = \{(a \oplus b) \odot c^* \mid b, c \in I\}$.
\item[(iii)] The ideal $(S]$ generated by $S$ is proper if and only if, for any $n < \omega$ and for any $a_1, \ldots, a_n \in S$, $a_1 \oplus \cdots \oplus a_n < 1$.
\end{enumerate}
\end{proposition}

In what follows, we shall always denote by $\widehat a$ and $\widehat X$, respectively, $\iota(a) \in [0,1]^{\Max A}$ and $\iota(X) \subseteq [0,1]^{\Max A}$, for $a \in A$ and $X \subseteq A$.

The class of semisimple MV-algebras form a full subcategory of $\MV$ that we shall denote by $\MVs$. As usual, for subsets $Z \subseteq Y$ of an ordered set $\la X \leq\ra$ we shall denote by $l_YZ$ (or simply $lZ$ when $Y = X$) the set of lower bounds of $Z$ in $Y$ and by $u_YZ$ (respectively: $uZ$) the set of all upper bounds of $Z$ in $Y$. We also recall that a subset $Y$ of $X$ is called a \emph{cut} if $Y = luY$. We set the following
\begin{definition}\label{ducuco}
Let $A$ be a semisimple MV-algebra. We say that a cut $X$ of $A$ is a \emph{limit cut} iff
\begin{equation}\label{distance}
%\begin{array}{r
 d(\widehat X,\widehat{uX}) = \bigwedge \{d(\widehat a, \widehat b) \mid b \in uX, a \in X\} = \bigwedge \{\widehat b \ominus \widehat a \mid b \in uX, a \in X\} = 0.
%\end{array}
\end{equation}

We shall say that $A$ is \emph{limit cut complete} (\emph{lcc} for short) if, for any limit cut $X$ of $A$, there exists in $A$ the supremum of $X$ or, equivalently, the supremum of $\widehat X$ in $[0,1]^{\Max A}$ belongs to $\widehat A$. 
\end{definition}

\begin{proposition}\label{lc}
Let $A$ be a semisimple MV-algebra. Then a cut $X$ of $A$ is a limit cut if and only if there exists a cut $Y$ of $A$ such that, in $[0,1]^{\Max A}$, $\bigvee \widehat X = \bigwedge \widehat Y^*$, where $Y^* = \{y^* \mid y \in Y\}$. Moreover, $Y$ is a limit cut too.
\end{proposition}
\begin{proof}
Let $X$ be a limit cut of $A$ and set $Y = (uX)^*$. From $x \leq y$ iff $x^* \geq y^*$ readily follows that $a \in uY$ iff $a^* \in luX = X$, whence $uY = X^*$. Analogously $a \in luY$ iff $a^* \in uX$. Therefore $luY = (uX)^* = Y$, i.e., $Y$ is a cut. Now, since $x \ominus y = 0$ iff $x \leq y$ in any MV-algebra, from $d(\widehat X, \widehat{uX}) = 0$, we get $\bigvee \widehat X = \bigwedge \widehat{uX} = \bigwedge \widehat Y^*$. Moreover, from $y^* \ominus x^* = y^* \odot x = x \ominus y$, we have that $d\left(\widehat Y,\widehat{uY}\right) = d\left(\widehat{(uX)}^*, \widehat X^*\right) = d(\widehat X, \widehat{uX}) = 0$, and therefore $Y$ is a limit cut. 

Conversely, let $X$ and $Y$ be cuts such that $\bigvee \widehat X = \bigwedge \widehat Y^*$, so in particular $d(\widehat X,\widehat Y) = 0$. Then $Y^* \subseteq uX$, whence $d(\widehat X,\widehat{uX}) \leq d(\widehat X,\widehat Y^*) = 0$, and $X$ is a limit cut. The fact that also $Y$ is a limit cut is an immediate consequence of the mutual roles of $X$ and $Y$ in this part of the proof.
\end{proof}

\begin{corollary}\label{lccch}
A semisimple MV-algebra $A$ is lcc if and only if, for all $X, Y \subseteq A$ and $\alpha \in [0,1]^{\Max A}$, $\alpha = \bigvee \widehat X = \bigwedge \widehat Y$ implies $\alpha \in \widehat A$.
\end{corollary}
\begin{proof}
Follows immediately from Proposition \ref{lc} by observing that, for any subset $X$ of $A$, $\bigvee \widehat X = \bigvee \widehat{luX}$. Then, if $\alpha = \bigvee \widehat X = \bigwedge \widehat Y$, $luX$ and $lu(Y^*)$ form a pair of limit cuts as in Proposition \ref{lc}.
\end{proof}

We wish to underline that the distance $d(\widehat X,\widehat{uX})$ considered in (\ref{distance}) do not necessarily coincide with $\iota(d(X,uX))$, as the following example shows.
\begin{exm}
Let $B$ the finite-cofinite Boolean algebra on the natural numbers. Let $\mathbb E$ be the set of even numbers and consider the set $X$ of all finite subsets of $\mathbb E$ and the set $Y$ of all cofinite subsets of $\N$ which include $\mathbb E$. Then it is self-evident that $X$ and $Y^*$ are cuts in $B$, $Y = uX$, and $d(X,Y)=0$ in $B$. However, by the Boolean Prime Ideal Theorem, we know that there exists a maximal ideal $M$ of $B$ which separates $X$ and $Y$, i.e. such that $X \subset M$ and $Y \cap M = \varnothing$. It follows that $d(\widehat X, \widehat Y) \neq 0$.
\end{exm}

The concept of limit cut complete MV-algebra arose naturally while the author was attempting to extend Stone duality to MV-algebras. Indeed, limit cut completeness is a distinctive feature of Boolean algebras among semisimple MV-algebras; in other words, all Boolean algebras are limit cut complete, while not all semisimple MV-algebras are. This circumstance shall appear clearer to the reader in the proof of the Duality Theorem. So, after all, the definition of limit cut complete MV-algebras is somehow ad hoc but, on the other hand, it turns out that the class of limit cut complete MV-algebras can play an important role for the theory of MV-algebras, as shown by the results of Section \ref{lccsec} and, in particular, by the fact that it is a reflective subcategory of $\MV$ and a completion subcategory of $\MVs$.

We shall try to describe as much as possible the class (in fact, the full subcategory of $\MVs$) $\MVlcc$ of limit cut complete MV-algebras in Section \ref{lccsec}. The rest of the present section is devoted to the extension of Stone Duality to MV-algebras and MV-topologies.

Let us now consider an MV-algebra $A$. By Theorem \ref{belrepr} and the comments following it, up to an isomorphism, $A' = A/Rad A$ is a subalgebra of $[0,1]^{\Max A}$. Therefore, $A'$ is a covering of $\Max A$ and, since it is an MV-subalgebra of $[0,1]^{\Max A}$, it is closed under $\oplus$, $\odot$ and $\wedge$. Then, by Proposition \ref{cov}, it is a base for an MV-topology $\O_A$ on $\Max A$. Conversely, given an MV-topological space $\tau = \la X, \O \ra$, the set $\clop\tau = \O \cap \Xi$ of the \emph{clopen} subsets of $X$, i.e. the fuzzy subsets of $X$ that are both open and closed, is a semisimple MV-algebra. Indeed $\0, \1 \in \clop\tau$ and, obviously, $\clop\tau$ is closed under $\oplus$ and $^*$; $\clop\tau$ is semisimple as an obvious consequence of being a subalgebra of $[0,1]^X$.

Let $\TMV$ be the category whose objects are MV-topological spaces and morphisms are MV-continuous functions between them. Moreover, we shall denote by $\SMV$ the full subcategory of $\TMV$ whose objects are \emph{Stone MV-spaces}, i.e., compact, separated MV-topological spaces having a base of clopen sets (\emph{zero-dimensional}).

In the proof of the following results we shall often identify any semisimple MV-algebra $A$ with its isomorphic image included in $[0,1]^{\Max A}$; so any element $a$ of a semisimple MV-algebra will be identified with the fuzzy set $\widehat a$. The reader may refer to \cite{bel,cha1,cha2,mvbook} for further details.

Let us now consider the following class functions:
\begin{equation}\label{func}
\begin{array}{cccccccccc}
\clop: & \tau & \in & \TMV & \ & \lmapsto & \ & \clop\tau & \in & \MV \\
\Max: & A & \in & \MV & & \lmapsto & & \la \Max A, \O_A\ra & \in & \TMV.
\end{array}
\end{equation}
Moreover, we set the following:
\begin{itemize}
\item for any two MV-topological spaces $\tau$ and $\tau'$, and for any continuous function $f: \tau \to \tau'$,
$$\clop f(\alpha) = f\fcou(\alpha), \textrm{ for all } \alpha \in \clop\tau';$$
\item for any two  MV-algebras $A$ and $B$, and for any MV-algebra homomorphism $h: A \to B$,
$$\Max h(N) = h\cou[N], \textrm{ for all } N \in \Max B.$$
\end{itemize}

\begin{lemma}
With the above notations, $\clop$ and $\Max$ are two contravariant functors.
\end{lemma}
\begin{proof}
%We will now show that the above mappings can be made into two contravariant functors. In order to do that, we need define, for any MV-algebra morphism $f: A \to B$, a continuous map $\Max f: \Max B \to \Max A$ and, for any continuous map $f: \tau \to \tau'$ between two arbitrary MV-spaces, an MV-algebra homomorphism $\clop f: \clop\tau' \to \clop\tau$. 

Let $\tau = \la X, \O_X \ra$ and $\tau' = \la Y, \O_Y \ra$ be two MV-topologies, and let $f: X \lto Y$ be a continuous map between them. As we already remarked $f\fcou: [0,1]^Y \lto [0,1]^X$ is a homomorphism of MV-algebras. On the other hand, by Definition \ref{cont}, $f\fcou[\O_Y] \subseteq \O_X$ and, as we observed right after the same definition, $f\fcou[\Xi_Y] \subseteq \Xi_X$; therefore $f\fcou[\clop\tau'] \subseteq \clop\tau$. Hence, for all $f \in \hom_{\TMV}(\tau,\tau')$, $\clop f$ is an MV-algebra homomorphism from $\clop \tau'$ to $\clop \tau$, i.e., a morphism from $\clop\tau$ to $\clop\tau'$ in $(\MVs)^{\op}$. The fact that $\clop$ preserves composition and identities is absolutely trivial.

Let now $A$ and $B$ be two MV-algebras and $h: A \lto B$ an MV-algebra homomorphism. It is known that the preimage of a maximal ideal under an MV-algebra homomorphism is a maximal ideal; then it is well-defined the map $\Max h: N \in \Max B \lmapsto h\cou[N] \in \Max A$. The function $\Max h$, on its turn, defines an MV-algebra homomorphism $(\Max h)\fcou: \alpha \in [0,1]^{\Max A} \lmapsto \alpha \circ \Max h \in [0,1]^{\Max B}$. Let us prove that $(\Max h)\fcou[A'] \subseteq \O_B$. 

So let $N$ be an arbitrary maximal ideal of $B$ and $M = \Max h(N)$. We have
$$\textrm{$(\Max h)\fcou(\widehat a)(N) = (\widehat a \circ \Max h)(N) = \widehat a(M)$, for all $a \in A$}.$$
The map $h': a/M \in A/M \lto h(a)/N \in B/N$ is well-defined since
$$\begin{array}{l}
a/M = a'/M \quad \To \quad (a \odot a'^*) \oplus (a' \odot a^*) \in M \quad \To \\
(h(a) \odot h(a')^*) \oplus (h(a') \odot h(a)^*) \in N \ \To \ h(a)/N = h(a')/N;
\end{array}$$
moreover it can be proved in a similar way that $h(a)/N = h(a')/N$ implies $a/M = a'/M$, that is, $h'$ is injective. Now, if we look at $A/M$ and $B/N$ as subalgebras of $[0,1]$, we get that the fuzzy set $\widehat{h(a)}$ takes, in any given $N \in \Max B$, precisely the same value taken by the fuzzy set $\widehat a$ in $M = \Max h(N)$. In other words, the fuzzy set $(\Max h)\fcou(\widehat a)$ is in $B'$, for all $a \in A$. It follows that $(\Max h)\fcou[A] \subseteq \O_B$ and therefore, by Lemma \ref{bascont}, $\Max h$ is a MV-continuous function from $\la \Max B, \O_B\ra$ to $\la \Max A, \O_A\ra$, i.e., it is a morphism from $\la \Max A, \O_A\ra$ to $\la \Max B, \O_B\ra$ in $\TMV^{\op}$. Again, it is immediate to see that $\Max$ is well-behaved w.r.t. composition and identity morphisms. 
\end{proof}

We recall (see \cite{mom,mvbook,mubook}) that an MV-algebra $A$ is called \emph{hyper-Archimedean} if all of its elements are Archimedean, namely, if all $a \in A$ satisfies the following equivalent conditions:
\begin{enumerate}[(a)]
\item there exists a positive integer $n$ such that $na \in \B(A)$;
\item there exists a positive integer $n$ such that $a^* \vee na = 1$;
\item there exists a positive integer $n$ such that $na = (n + 1)a$.
\end{enumerate}
It is well-known that every hyper-Archimedean MV-algebra is semisimple while the converse is not true. Moreover, an MV-algebra $A$ is hyper-Archimedean if and only if it is isomorphic to a Boolean product of subalgebras of $[0,1]$ (the reader may refer to \cite[Section 6.5]{mvbook} for more details). As to the relationship between hyper-Archimedean and lcc MV-algebras, it must be noticed that none of the two classes is included in the other one. Indeed, for example, $[0,1] \cap \Q$ is hyper-Archimedean and not lcc, while any algebra of type $[0,1]^X$, with $X$ infinite, is easily seen to be lcc and not hyper-Archimedean. However, the two classes have a non-trivial intersection which includes the whole class of liminary MV-algebras (see Definition \ref{liminary}), as shown in the last section of the paper.

\begin{theorem}[Duality theorem]\label{dual}
$\clop$ and $\Max$ form a duality between $\MVlcc$ and $\SMV$.
\end{theorem}
\begin{proof}
It is immediate to verify that both the functors, restricted to $\SMV$ and $\MVlcc$ respectively, are faithful. We shall prove that
$$\Max\clop \tau \cong_{\TMV} \tau \quad \textrm{ and } \quad \clop\Max A \cong_{\MV} A,$$
for all $\tau \in \SMV$ and for all $A \in \MVlcc$. The assertion will therefore follow from the fact that such isomorphisms, together with faithfulness, yield two natural isomorphisms between the two compositions $\Max\clop$ and $\clop\Max$ and, respectively, $\id_{\SMV}$ and $\id_{\MVlcc}$.

First, let us prove that $\Max A \in \SMV$ for any semisimple MV-algebra $A$ and that $\clop\Max A \cong A$ if $A \in \MVlcc$.

$\la \Max A, \O_A\ra$ is zero-dimensional by definition. $\clop\Max A$ is obviously semisimple, and every element of $\clop\Max A$ can be obtained as both a join and a meet of elements of $\widehat A$. Therefore, if $A \in \MVlcc$, by Proposition \ref{lc} and Corollary \ref{lccch}, $A \cong \clop\Max A$.  Now we need to prove only that $\Max A$ is compact and Hausdorff. Let $\Gamma$ be an open covering of $\Max A$ and assume, by contradiction, that it does not contain any additive covering. By Proposition \ref{properid}, $(\Gamma]$ is a proper ideal of $A$ and, therefore, it is contained in some $M \in \Max A$; but this implies that, for any $a \in \Gamma$, $a(M) = 0$, i.e. $\Gamma$ is not  a covering of $\Max A$, which is absurd.

In order to prove separation, let us consider $M \neq N \in \Max A$ and let $a \in M \setminus N$; we have $\widehat a(M) = 0$ and $\widehat a(N) \neq 0$. So, since $[0,1]$ is hyper-Archimedean, there exists $k < \omega$ such that $k \widehat a(N) = 1$. % and, by Proposition \ref{maxch}, there exists $n < \omega$ such that $(\widehat a^*)^n \in N$, that is, $(\widehat a^*)^n(N) = 0$. If $k = \max\{m,n\}$, 
Then we have $k \widehat a(N) = 1$ and $\widehat a^*(M) = \widehat a(M)^* = 1$, which implies $(\widehat a^*)^k(M) = 1$. %, and $\widehat a(M) = 0$ implies $k \widehat a(M) = 0$. Therefore $(\widehat a^*)^k$ and $k\widehat a$ are two open sets satisfying (i) and (ii) of Definition \ref{t2}. 
Moreover, $(\widehat a^*)^k \odot k\widehat a = (k\widehat a)^* \odot k\widehat a = \0$; then, by Remark \ref{t2r}, $\la \Max A, \O_A\ra$ is a Stone MV-space.

Now let us prove that $\tau$ and $\Max\clop\tau$ are homeomorphic for any Stone MV-space $\tau$. Let $\tau = \la X, \O \ra$ be a Stone MV-space and, for each $x \in X$, let $f(x) = \{o \in \clop\tau \mid o(x) = 0\}$. It is self-evident that $f(x)$ is a proper ideal of the algebra $A = \clop\tau$ for all $x \in X$. For any fixed $x$ and for each $o \in A$, $o \notin f(x)$ implies $o(x) > 0$ and, therefore, $o^*(x) < 1$. Then there exists $n < \omega$ such that $(o^*)^n(x) = 0$, i.e. $(o^*)^n \in f(x)$, and Proposition \ref{maxch} ensures us that $f(x)$ is a maximal ideal.

Now we must prove that the map $f: X \lto \Max A$ is a homeomorphism of MV-spaces. First, let $x \neq y \in X$; since $\tau$ is Hausdorff, there exist $o_x, o_y \in \O$ that satisfy Definition \ref{t2}, and each of these open sets is the join of a set of clopens because $\tau$ is zero-dimensional. By Lemma \ref{xclos}, $\{x\}$ and $\{y\}$ are closed, whence, by Lemma \ref{closcomp}, they are compact; then there exist two finite families of such sets --- say $\{o_{xi}\}_{i=1}^n$ and $\{o_{yj}\}_{j=1}^m$ --- which are additive open coverings of $\{x\}$ and $\{y\}$ respectively, and are such that $(o_{x1} \oplus \cdots \oplus o_{xn})(y)= 0 = (o_{y1} \oplus \cdots \oplus o_{ym})(x)$. Moreover, $o_{x1} \oplus \cdots \oplus o_{xn}$ and $o_{y1} \oplus \cdots \oplus o_{ym}$ are both clopen, hence the former belongs to $f(y)$ and the latter to $f(x)$. It follows $f(x) \neq f(y)$, namely, $f$ is injective.

In order to prove that $f$ is onto, let $M \in \Max A$ and assume, by contradiction, that $M$ is not the image under $f$ of any element of $X$, that is, for all $x \in X$ there exists $o \in M$ such that $o(x) > 0$. Then, for each $x \in X$, there exist $o \in M$ and $m < \omega$ such that $mo(x) = 1$, and $mo \in M$ because $M$ is an ideal. So let, for each $x \in X$, $o_x$ be an element of $M$ whose value in $x$ is $1$; the family $\{o_x\}_{x \in X}$ is an open covering of $X$ whence, by the compactness of $\tau$, it contains an additive covering $\{o_i\}_{i=1}^n$. It follows that $\1 = o_1 \oplus \cdots \oplus o_n \in M$ which contradicts the hypothesis that $M$ is a proper ideal. Such a contradiction follows from the assumption that for all $x \in X$ there exists $o \in M$ such that $o(x) > 0$; hence there exists $x \in X$ such that $o(x) = 0$ for all $o \in M$, i.e., such that $M = f(x)$, and $f$ is onto.

We need to prove that both $f$ and $f^{-1}$ are continuous. To this purpose, we first observe that, for all $x \in X$ and $o \in \clop\tau$, $o/f(x)$ is a real number in $[0,1]$ and coincide with the membership value $o(x)$ of the point $x$ to the clopen $o$. Indeed, by Proposition \ref{classch}, $o/f(x) = \{(o \oplus p) \odot q^* \mid p, q \in f(x)\}$ and, on the other hand, $((o \oplus p) \odot q^*)(x) = (o(x) \oplus 0) \odot 1 = o(x)$ for all $p, q \in f(x)$. Therefore, $\clop\tau/f(x) = \{o(x) \mid o \in \clop\tau\}$ and $\pi_{f(x)}: o \in \clop\tau \mapsto o(x) \in \clop\tau/f(x) \subseteq [0,1]$.

Now, any clopen $o$ of $\tau$ can be identified (see the proof of Theorem \ref{belrepr}) with a clopen $\widehat o$ of $\Max\clop\tau$ in a unique way: $\widehat o(M) = \iota_M(\pi_M(o)) = \iota_{f(x)}(\pi_{f(x)}(o)) = \iota_{f(x)}(o(x))$, for all $M = f(x) \in \Max\clop\tau$, and $\iota_{f(x)}$ is simply the inclusion map of $\clop\tau/f(x)$ in $[0,1]$. Therefore, for any basic clopen $\widehat o$ of $\Max\clop\tau$, and for each $x \in X$, $f\fcou(\widehat o)(x) = (\widehat o \circ f)(x) = \widehat o(f(x)) = o(x)$, with $o \in \clop\tau$. It follows that the fuzzy preimage, under $f$, of any basic open set of $\Max\clop\tau$ is open in $\tau$, that is, $f$ is continuous. Analogously, for each $M = f(x) \in \Max\clop\tau$, $(f^{-1})\fcou(o)(M) = (o \circ f^{-1})(f(x)) = o(x) = \widehat o(M)$, and $f^{-1}$ is continuous as well. We can conclude that $\tau$ and $\Max\clop\tau$ are homeomorphic spaces.

The proof is complete.
\end{proof}

\begin{corollary}\label{restr0}
The restriction of the above duality to Boolean algebras and crisp topologies coincide with the classical Stone Duality. 
\end{corollary}
\begin{proof}
This is a trivial consequence of how the functors are defined.
\end{proof}

\begin{theorem}\label{restr}
For any Stone MV-space $\tau$, its skeleton space is a Stone space and its image under $\clop$ is precisely the Boolean center of $\clop\tau$.

Conversely, for any semisimple MV-algebra $A$, $\Max\B(A)$ coincide with the skeleton topology of $\Max A$.
\end{theorem}
\begin{proof}
The first part is trivial. For the second part, once observed that, by Proposition \ref{maxch}, $M \cap \B(A)$ is a maximal ideal of the Boolean algebra $\B(A)$, for all $M \in \Max A$, it suffices to apply Theorem \ref{dual}.
\end{proof}

It is immediate to verify that
$$\begin{array}{cccc}
\B: & A \in \MV & \lmapsto & \B(A) \in \BB \\
\sk: & \la X, \O\ra \in \TMV & \lmapsto & \la X, \B(\O) \ra \in \Top
\end{array}$$
define two functors, where the action of $\B$ on morphisms is simply the restriction of the MV-algebra homomorphism to the Boolean center of the domain, and $\sk f$ is $f$ itself, for any MV-continous map $f$. They are, in fact, the left-inverses of the inclusion functors. Then Theorem \ref{restr} (together with Corollary \ref{restr0}) can be reformulated as follows.

\begin{corollary}\label{restr2}
$\clop_\restr \circ \sk = \B \circ \clop$ and $\Max_\restr \circ \B = \sk \circ \Max$.
\end{corollary}

Then we have the following commutative diagram of functors, where horizontal arrows are equivalences and vertical ones are inclusions of full subcategories and their respective left-inverses.
\begin{equation}\label{diagfunc}
\xymatrix{
\MVlcc \ar@<.5ex>[rr]^\Max \ar@<-.5ex>[dd]_\B & & \SMV^{\op} \ar@<.5ex>[dd]^\sk \ar@<.5ex>[ll]^\clop \\
&&\\
\BB \ar@<-.5ex>[uu]_{\rotatebox[origin=c]{90}{$\subseteq$}} \ar@<-.5ex>[rr]_{\Max_\restr} & & \ST^{\op} \ar@<.5ex>[uu]^{\rotatebox[origin=c]{90}{$\subseteq$}} \ar@<-.5ex>[ll]_{\clop_\restr} \\
}
\end{equation}

\section{On limit cut complete MV-algebras}\label{lccsec}

In the present section we shall describe the class of limit cut complete MV-algebras, namely, the category which is dual to the one of Stone MV-spaces, and we will show that it is a completion subcategory of $\MVs$ and therefore a reflective subcategory of $\MV$. In order to do that, we will show that, for any semisimple MV-algebra $A$, the extension $\clop\Max A$ is the smallest lcc MV-algebra containing $A$, and it can be obtained also with an alternative construction within the class of MV-algebras itself.

We already presented a characterization of $\MVlcc$ (Theorem \ref{lccch}) which, however, may not be handy enough in many cases. As we shall see, a necessary condition for a semisimple MV-algebra to be lcc is that all of its quotients on the maximal ideals must be complete chains or, equivalently, either finite or isomorphic to $[0,1]$. Whether such a condition is sufficient too, is still an open problem, as the author was able neither to prove nor to disprove it so far. It will be stated as a conjecture at the end of the section.

We shall begin by showing two important properties of lcc MV-algebras.

\begin{theorem}\label{lccnec}
If $A$ is an lcc MV-algebra, then the following hold.
\begin{enumerate}
\item[(i)] For all $a \leq b \in A$, if $\iota([a,b])$ is dense in $[\widehat a, \widehat b]$, then $\iota([a,b]) = [\widehat a, \widehat b]$.
\item[(ii)] $A$ is a subdirect product of complete MV-chains, i.e., for all $M \in \Max A$, $A/M$ is isomorphic either to $[0,1]$ or to the finite chain $\textrm{\emph{\L}}_n$ for some $n < \omega$.
\end{enumerate}
\end{theorem}
\begin{proof}
\begin{enumerate}
\item[(i)] If $A$ is lcc and $\iota([a,b])$ is dense in $[\widehat a, \widehat b]$, then every element of $[\widehat a, \widehat b]$ can be obtained as both a supremum and an infimum of elements of $\iota([a,b])$. Then, by Corollary \ref{lccch}, $[\widehat a, \widehat b] = \iota([a,b])$.
\item[(ii)] Let us consider the dual space $\la \Max A, \Omega_A\ra$ of $A$. Then, for all $M \in \Max A$, the subspace $\la \{M\}, \Omega_{A,M} \ra$ is clearly a Stone MV-space. Now, since the (continuous) inclusion map $\{M\} \to \Max A$ corresponds, by Theorem \ref{dual}, to the canonical projection $\pi_M: A \to A/M$, the quotient $A/M$ must be isomorphic to $\clop\{M\}$. On the other hand, it is easy to see that the only possible MV-topologies on a singleton are either $[0,1]$ or $\textrm{\L}_n$ for some $n < \omega$. So, if $A/M$ is infinite, it must be dense in $[0,1]$, by \cite[Proposition 3.5.3]{mvbook}, and $\{\bigvee X \mid X \subseteq A/M\} = [0,1]$, that is, $A/M \cong [0,1]$.
\end{enumerate}
\end{proof}

In the next three results we shall prove that liminary MV-algebras are lcc, and that they are dual to strongly compact Stone MV-spaces. First, we recall the following 
\begin{definition}\label{liminary}
An MV-algebra $A$ is called \emph{liminary} if all of its quotients on the prime ideals are finite \cite{cielmu}.
\end{definition}
It is immediate to see that all liminary MV-algebras are locally finite (namely, all of their finitely generated subalgebras are finite), hence hyper-Archimedean (see, for instance, \cite[Section 6.5]{mvbook}.
\begin{proposition}\label{lflem}
Let $A$ be a liminary MV-algebra. Then $\Max A$ is strongly compact.
\end{proposition}
\begin{proof}
Let $A$ be liminary and $\Gamma \subseteq \widehat A$ be a covering of $\Max A$ made of basic clopens. Since all the quotients of $A$ over maximal ideals are finite chains, it follows immediately that, for each $M \in \Max A$, there exists an element $\widehat a \in \Gamma$   such that $\widehat a(M) = 1$. So, for each maximal ideal $M$, let $\widehat a_M \in \Gamma$ be one of such clopens; clearly $\{\widehat{a_M}\}_{M \in \Max A}$ is a subcovering of $\Gamma$. Since $A$ is liminary, it is hyper-Archimedean, and therefore, for all $M \in \Max A$, there exists $n_M < \omega$ such that $\widehat{a_M^{n_M}}$ is Boolean. Then the family $\{\widehat{a_M^{n_M}}\}_{M \in \Max A}$ is, again, a covering of $\Max A$. Since $\Max A$ is compact, there exist $M_1, \ldots, M_k \in \Max A$ such that $\bigoplus_{i = 1}^k \widehat{a_{M_i}^{n_{M_i}}} = \1$. But the clopens of type $\widehat{a_M^{n_M}}$ are Boolean, and therefore we have
$$\1 = \bigoplus_{i = 1}^k \widehat{a_{M_i}^{n_{M_i}}} = \bigvee_{i = 1}^k \widehat{a_{M_i}^{n_{M_i}}} \leq \bigvee_{i = 1}^k \widehat{a_{M_i}}.$$
Hence $\{\widehat{a_{M_i}}\}_{i=1}^k$ is a finite subcovering of $\Gamma$ and $\Max A$ is strongly compact.
\end{proof}

\begin{proposition}\label{lflem2}
Liminary MV-algebras are limit cut complete.
\end{proposition}
\begin{proof}
Let $A$ be a liminary MV-algebra and let $\alpha \in [0,1]^{\Max A}$ be such that there exist $\widehat X, \widehat Y \subseteq \widehat A$ with $\alpha = \bigvee \widehat X = \bigwedge \widehat Y$. Then, as in Proposition \ref{lc}, $d(\widehat X, \widehat Y) = \bigwedge \{\widehat b \ominus \widehat a \mid b \in Y, a \in X\} = 0$. Therefore $\{\widehat b^* \oplus \widehat a \mid b \in Y, a \in X\}$ is a covering of $\Max A$ and, by Proposition \ref{lflem}, it contains a finite covering, i.e., there exist $a_1, \ldots, a_k \in X$ and $b_1, \ldots, b_k \in Y$ such that $\bigvee_{i=1}^k \widehat b_i^* \oplus \widehat a_i = \1$. On the other hand, $\bigvee_{i=1}^k \widehat b_i^* \oplus \widehat a_i \leq \left(\left(\bigwedge_{i=1}^k \widehat b_i\right)^*\right) \oplus \left(\bigvee_{i=1}^k \widehat a_i\right)$, with $\bigwedge_{i=1}^k \widehat b_i \in \widehat Y$ and $\bigvee_{i=1}^k \widehat a_i \in \widehat X$, whence $\alpha = \bigvee_{i=1}^k \widehat a_i = \widehat{\bigvee_{i=1}^k  a_i} \in \widehat A$. The thesis follows.
\end{proof}

\begin{proposition}
If $\tau = \la X, \Omega\ra$ is a strongly compact Stone MV-space, then $A = \clop \tau$ is a liminary MV-algebra.
\end{proposition}
\begin{proof}
First of all, it is immediate to see that $A$ is a subdirect product of finite chains. Indeed, combining by \ref{lccnec}(ii) with strong compactness, we get immediately that, for any $x \in X$, $\clop \{x\}$, which is the quotient of $A$ by a maximal ideal, is necessarily a finite chain, whence $A/M$ is a finite chain for all $M \in \Max A$.

Let us now consider an arbitrary element $\alpha$ of $A$. Then, for all $x \in \supp \alpha$, there exists $n_x < \omega$ such that $n_x\alpha(x) = 1$ and, for each $x \in \supp\alpha$, $(\alpha^*)^{n_x}(x) = 0$ and $(\alpha^*)^{n_x}(y) = 1$ for all $y \in X \setminus \supp\alpha$. So the family $\{n_x\alpha\}_{x \in \supp\alpha} \cup \{(\alpha^*)^{n_x}\}_{x \in \supp\alpha}$ is an open covering of $X$. Since $X$ is strongly compact, there exist $x_1, \ldots, x_k, y_1, \ldots, y_h \in \supp\alpha$ such that $\left(\bigvee_{i=1}^k n_{x_i}\alpha\right) \vee \left(\bigvee_{j=1}^h (\alpha^*)^{n_{y_j}}\right) = \1$. Now, if we set $n = \max\{n_{x_i}\}_{i=1}^k$, we obtain that $n\alpha(x) = \left\{\begin{array}{l} 1 \textrm{ if $x \in \supp\alpha$} \\  0 \textrm{ if $x \in X \setminus \supp\alpha$} \end{array}\right.$, namely, $n\alpha \in \B(A)$. It follows that any element of $A$ is Archimedean, whence $A$ is an hyper-Archimedean algebra whose quotients on the maximal ideals are all finite. Then $A$ is liminary.
\end{proof}

%We will now consider the case of \emph{weakly locally finite MV-algebras} \cite{cima}, i. e., those MV-algebras whose finitely generated subalgebras are finite direct products of subalgebras of $[0,1]$. 

Let now, for any semisimple MV-algebra $A$, $\cat{LCC}(A)$ be the following set
$$\cat{LCC}(A) = \{B \leq [0,1]^{\Max A} \mid B \in \MVlcc \textrm{ and } \widehat A \leq B\}.$$
$\cat{LCC}(A)$ is certainly not empty since it contains at least the Dedekind-MacNeille completion $A^{DM}$ of $A$ (see \cite[Section 5]{bgl}) and $[0,1]^{\Max A}$ (which may possibly coincide). Then we can set the following
\begin{definition}\label{lccclo}
For any semisimple MV-algebra $A$, the MV-algebra $A^{\lcc} = \bigcap \cat{LCC}(A)$ will be called the \emph{limit cut completion} (\emph{lc-completion}, for short) of $A$.
\end{definition}

Since, for any limit cut $X$ of $A$, $\widehat X$ is a limit cut of $\widehat A$, and $\bigvee \widehat X \in B$ for all $B \in \cat{LCC}(A)$, we get immediately that $A^{\lcc} \in \MVlcc$ for any $A \in \MVs$.

\begin{theorem}\label{lccfunc}
The mapping
$$(~~)^{\lcc}: A \in \MVs \mapsto A^{\lcc} \in \MVlcc,$$
with $f^{\lcc} = \clop\Max f$ for any morphism $f$ of $\MVs$, defines a categorical completion, namely, a faithful reflection.
\end{theorem}
\begin{proof}
From the results of Section \ref{duality} it readily follows that, for all $A \in \MVs$, $A^{\lcc} = \clop\Max A$. Then, by Theorem \ref{dual}, $(~~)^{\lcc}$ is left adjoint to the inclusion functor, whose faithfulness is obvious.
\end{proof}

So, the category $\MVlcc$ is a completion of $\MVs$, i.e., a reflective subcategory whose reflector is faithful. Therefore, in particular, we have the following immediate consequence.
\begin{corollary}\label{univ}
For any semisimple MV-algebra $A$, the algebra $A^{\lcc}$ has the following universal property: for any lcc MV-algebra $B$ and for any homomorphism $f: A \to B$ there exists a unique homomorphism $f^{\lcc}: A^{\lcc} \to B$ such that $f^{\lcc}_{\restr A} = f$.

Equivalently, for any semisimple MV-algebra $B$ and for any homomorphism $f: A \to B$ there exists a unique homomorphism $f^{\lcc}: A^{\lcc} \to B^{\lcc}$ such that $f^{\lcc}(a) = f(a)$ for all $a \in A$.
\end{corollary}

Now taking into account that $\MVs$ is a reflective subcategory of $\MV$, we obtain also the following
\begin{corollary}\label{refl}
$\MVlcc$ is a reflective subcategory of $\MV$.
\end{corollary}

\begin{definition}
We shall say that an MV-algebra $A$ is \emph{subdirect factor complete} iff it is a subdirect product of complete chains. In what follows, we shall denote by $\MVsfc$ the full subcategory of $\MVs$ whose objects are subdirect factor complete algebras. 
\end{definition}

By Theorem \ref{lccnec}, the class $\MVlcc$ is contained in $\MVsfc$. We set the following
\begin{conjecture}
$\MVlcc = \MVsfc$, namely, every subdirect factor complete MV-algebra is limit cut complete and, therefore, is the dual algebra of a Stone MV-space. 
\end{conjecture}

%\begin{corollary}\label{hyper}
%Stone MV-spaces which are strongly separated are dual to dense-complete hyper-Archimedean MV-algebras.
%\end{corollary}
%\begin{proof}
%With reference to last part of the proof of Claim 2 in Theorem \ref{dual}, it suffices to observe that, according to Definition \ref{arch}, the natural number $k$ can be chosen in such a way that $(a^*)^k$ and $ka$ are Boolean. Then we obtain $(a^*)^k \odot ka = (a^*)^k \wedge ka = \0$, whence $\la \Max A, \O_A\ra$ is strongly separated.
%\end{proof}

\section{Combining dualities}
\label{sec:comb}

The results presented in the previous sections provide, in our opinion, a strong motivation for the development of a more general and comprehensive theory of MV-topologies. Indeed, besides proving once more that MV-algebras are the most natural generalization of Boolean algebras, these results can be combined with the wide variety of equivalences involving categories of MV-algebras, thus giving new equivalences and, therefore, new tools.

For example, we can apply the well-known and celebrated categorical equivalences \cite{mun} between MV-algebras and lattice-ordered Abelian groups with a strong order unit (Abelian $u\ell$-groups, for short), thus obtaining a duality between a suitable category of Archimedean Abelian $u\ell$-groups and Stone MV-spaces.

It is worth noticing that, on their turn, Archimedean Abelian $u\ell$-groups are, up to isomorphisms, subgroups of the $u\ell$-group of bounded functions from a set $X$ to $\R$, with pointwise operations and the $1$-constant map as order unit. Then the restriction of such functors yields a duality between Stone spaces and $u\ell$-groups which are, up to isomorphisms, subgroups of the $u\ell$-group of bounded functions from a set $X$ to $\Z$.

Another example can be given with reference to \cite{dnl}, where the authors proved, for each $n > 1$, a categorical equivalence between the variety of MV-algebras generated by the $(n+1)$-element chain $\textrm{\L}_{n+1} = \{i/n)\}_{i=0}^n$ and the category whose objects are pairs $(B, R)$, where $B$ is a Boolean algebra and $R$ is an $n$-ary relation on $B$ satisfying certain conditions, and a morphism $f: (B,R) \lto (B',R')$ is a Boolean algebra homomorphism such that $(a_0, \ldots, a_{n-1}) \in R$ implies $(f(a_0), \ldots, f(a_{n-1})) \in R'$.  So, since all the MV-algebras in such varieties are liminary, we can combine our results with the ones in \cite{dnl} thus obtaining a ctegorical equivalence between $n$-valued Stone MV-spaces and a suitable category of classical Stone spaces with additional conditions.

Let us describe the situation in full details.

\begin{definition}\label{bbn}
We define the category $\BB_n$ as follows.
\begin{itemize}
\item An object of $\BB_n$ is a pair $B_n = \la B, (J_i)_{i=1}^{n-1}\ra$ where $B$ is a Boolean algebra and $(J_i)_{i=1}^{n-1}$ is a sequence of $n-1$ ideals of $B$ such that
\begin{enumerate}[(i)]
\item $J_i = J_{n-i}$ for all $i = 1, \ldots, n-1$, and
\item $J_h \cap J_{i-h} \subseteq J_i$, for all $i = 2, \ldots, n-1$ and $h = 1, \ldots, i-1$.
\end{enumerate}
\item For any two objects $B_n = \la B, (J_i)_{i=1}^{n-1}\ra$ and $B_n' = \la B', (J_i')_{i=1}^{n-1}\ra$ a morphism $f: B_n \lto B_n'$ is a Boolean algebra homomorphism from $B$ to $B'$ such that $f[J_i] \subseteq J_i'$ for all $i = 1, \ldots, n-1$. 
\end{itemize}
\end{definition}

In \cite{dnl} the authors defined the category $\cat BR_n$, for each $n \in \omega$. The objects of this category are pairs $\la B, R_n\ra$ where $B$ is a Boolean algebra and $R_n$ is an $n$-ary relation on $B$ such that
\begin{itemize}
\item $(a_0, \ldots, a_{n-1}) \in R_n$ implies $a_0 \geq a_1 \geq \cdots \geq a_{n-1}$,
\item $(a_0, \ldots, a_{n-1}) \in R_n$ implies $(a_{n-1}^*, \ldots, a_0^*) \in R_n$,
\item $\ul a = (a, a, \ldots, a) \in R_n$ for all $a \in B$, and
\item $(a_0, \ldots, a_{n-1}), (b_0, \ldots, b_{n-1}) \in R_n$ implies
$$\left(a_i \vee b_i \vee \bigvee_{h+k=i-1} (a_h \wedge b_k)\right)_{i=0}^{n-1} \in R_n.$$
\end{itemize}
Given two objects $\la B, R_n\ra$ and $\la B', R_n'\ra$ in $\cat BR_n$, a morphism $f$ between them is a Boolean algebra homomorphism $f: B \lto B'$ such that $(a_0, \ldots, a_{n-1}) \in R_n$ implies $(f(a_0), \ldots, f(a_{n-1})) \in R_n'$. 

\begin{lemma}\label{relid}
The categories $\BB_n$ and $\cat BR_n$ are isomorphic.
\end{lemma}
\begin{proof}
By \cite[Proposition 24]{dnl}, for any $\la B, R_n\ra \in \cat BR_n$,
\begin{equation}\label{jirn}
J_i(R_n) = \{a_{i-1} \wedge a_i^* \mid (a_0, \ldots, a_{n-1}) \in R_n\}
\end{equation}
is an ideal of $B$ for all $i = 1, \ldots, n-1$ and the sequence $(J_i(R_n))_{i=1}^{n-1}$ satisfies the conditions of Definition \ref{bbn}. Conversely, by \cite[Proposition 25]{dnl}, if $\la B, (J_i)_{i=1}^{n-1}\ra$ is an object of $\BB_n$, the set
\begin{equation}\label{rj}
R_J = \{(a_0, \ldots, a_{n-1}) \mid a_{i-1} \geq a_i \textrm{ and } a_{i-1} \wedge a_i^* \in J_i, \textrm{ for } i = 1, \ldots, n-1\}
\end{equation}
is an $n$-ary relation on $B$ that makes $\la B, R_J\ra$ an object of $\cat BR_n$. Moreover, by \cite[Lemma 26]{dnl} $R_n = R_{J(R_n)}$ for any $\la B, R_n \ra \in \cat BR_n$.

Reciprocally, let $\la B, (J_i)_{i=1}^{n-1}\ra \in \BB_n$; we shall prove that $J_i(R_J) = J_i$ for all $i$. The inclusion $J_i(R_J) \subseteq J_i$ follows immediately from (\ref{jirn}) and (\ref{rj}). On the other hand, for any index $i$ and for any element $a \in J_i$, the $n$-tuple
$$(1, 1, \ldots, \underbrace{1}_{i-1},\underbrace{a^*}_{i}, a^*, \ldots, a^*)$$
is clearly in $R_J$. Therefore $a \in J_i(R_J)$ whence $J_i(R_J) = J_i$.

Last, we need to prove that these mappings are functorial. So let us consider a morphism $f: \la B, R_n\ra \lto \la B', R_n'\ra$ in $\cat BR_n$. Since $f$ is also a Boolean algebra homomorphism between $B$ and $B'$ such that $(a_0, \ldots, a_{n-1}) \in R_n$ implies $(f(a_0), \ldots, f(a_{n-1})) \in R_n'$, we have that $f(a_{i-1} \wedge a^*_i) = f(a_{i-1}) \wedge f(a_i)^* \in J_i(R_n')$ for all $i = 1, \ldots, n-1$ and for all $(a_0, \ldots, a_{n-1}) \in R_n$. So $f[J_i(R_n)] \subseteq J_i(R_n')$ for all $i$, and $f$ is a morphism in $\BB_n$ from $\la B, (J_i(R_n))_{i=1}^{n-1}\ra$ to $\la B', (J_i(R_n'))_{i=1}^{n-1}\ra$. The proof of the other implication is completely analogous, hence $\BB_n$ and $\cat BR_n$ are isomorphic categories.
\end{proof}

%Thanks to Lemma \ref{relid} we can use indifferently 

\begin{definition}\label{topn}
We define the category $\Top_n$ as follows.
\begin{itemize}
\item An object of $\Top_n$ is a pair $\tau_n = \la \la X, \Omega \ra, (o_i)_{i=1}^{n-1}\ra$ where $\la X, \Omega \ra$ is a topological space and $(o_i)_{i=1}^{n-1}$ is a sequence of $n-1$ open subsets of $X$ such that
\begin{enumerate}[(i)]
\item $o_i = o_{n-i}$ for all $i = 1, \ldots, n-1$, and
\item $o_h \cap o_{i-h} \subseteq o_i$, for all $i = 2, 3, \ldots, n-1$ and $h = 1, \ldots, i-1$.
\end{enumerate}
\item For any two objects $\tau_n = \la \la X, \Omega \ra, (o_i)_{i=1}^{n-1}\ra$ and $\tau_n' = \la \la X', \Omega' \ra, (o_i')_{i=1}^{n-1}\ra$ a morphism $f: \tau_n \lto \tau_n'$ is a continuous map from $X$ to $X'$ such that $f\cou[o_i'] \subseteq o_i$ for all $i = 1, \ldots, n-1$. 
\end{itemize}
We shall denote by $\ST_n$ the full subcategory of $\Top_n$ whose objects have a Stone space as the underlying topology.
\end{definition}

\begin{theorem}\label{stonen}
The categories $\BB_n$ and $\ST_n$ are dual to each other.
\end{theorem}
\begin{proof}
By \cite[Theorem 7.25]{kop}, ideals of a Boolean algebra and open sets of its dual Stone space are dual to each other. Reformulating that result with our notations, given a Boolean algebra $B$ and its dual Stone space $\la \Max B, \Omega\ra$, the duality between the two structures defines an order isomorphism between the posets $\la \Id(B), \subseteq \ra$ and $\la\Omega, \subseteq\ra$. Indeed, for any ideal $I$ of $B$ the set $o_I = \bigvee_{a \in I} \widehat a$ is open in the dual space $\Max B$ of $B$. Conversely, for any Stone space $\la X, \Omega\ra$ and any open set $o$ of $X$, the set $I_o = \{a \in \clop X \mid a \leq o\}$ is an ideal of the Boolean algebra $\clop X$. Moreover, these two maps are order-preserving, bijective and inverses of each other, namely, for all $o, o' \in \Omega$ and $I,I' \in \Id(B)$ the following hold:
\begin{itemize}
\item $o \leq o'$ implies $I_o \subseteq I_{o'}$,
\item $I \subseteq I'$ implies $o_I \leq o_{I'}$,
\item $o_{I_o} = o$, and
\item $I_{o_I} = I$.
\end{itemize}
So we can define
\begin{equation}\label{ultn}
\Max_n: B_n = \la B, (J_i)_{i=1}^{n-1} \ra \in \BB_n \lmapsto \la \la \Max B, \Omega\ra, (o_{J_i})_{i=1}^{n-1}\ra \in \ST_n
\end{equation}
\begin{equation}\label{clopn}
\clop_n: \tau_n = \la \la X, \Omega \ra, (o_i)_{i=1}^{n-1} \ra \in \ST_n \lmapsto \la \clop X, (I_{o_i})_{i=1}^{n-1}\ra \in \BB_n.
\end{equation}
The previous discussion ensures that $\clop_n\Max_n B_n \cong B_n$ and $\Max_n\clop_n \tau_n \cong \tau_n$ in the respective categories. Let us prove that the two mappings really define two contravariant functors.

Now, given a morphism $f: B_n \lto B_n'$ in $\BB_n$, we already know that $f$ is in particular a Boolean algebra homomorphism from $B$ to $B'$; therefore $\Max_n f$ is the continuous map $ M' \in \Max B' \lmapsto f\cou[M'] \in \Max B$.  Therefore we just need to prove that $(\Max_n f)\cou[o_{J_i}] \leq o_{J_i'}'$.

For any index $i$ and for all $M' \in \Max B'$, we have
$$(\Max_n f)\cou[o_{J_i}](M') = (o_{J_i} \circ (\Max_n f))(M') = \left\{\begin{array}{ll} 0 & \textrm{ if $J_i \cap (\Max_n f)(M') \neq \varnothing$} \\ 1 & \textrm{ if $J_i \cap (\Max_n f)(M') = \varnothing$} \end{array}\right.$$
But, for all $a \in B$, $a \in (\Max_n f)(M')$ iff $a \in f\cou(M')$ iff $f(a) \in M'$, hence
$$(\Max_n f)\cou[o_{J_i}] = o_{J_i} \circ (\Max_n f) = \bigvee_{a \in J_i} \widehat{f(a)} \leq \bigvee_{a' \in J_i'} \widehat{a'} = o'_{J_i'}$$

Analogously, for any morphism $g: \tau_n \lto \tau_n'$, $\clop_n g$ is the Boolean algebra homomorphism $o' \in \clop X' \lmapsto g\cou[o'] \in \clop X$, and what we need to prove is that $\clop_n g [I_{o_i'}'] \subseteq I_{o_i}$ for all $i = 1, \ldots, n-1$. For all $i = 1, \ldots, n-1$, $\clop_n g [I_{o_i'}'] = \{g\cou[a'] \mid a' \in \clop X' \textrm{ and } a' \leq o_i'\}$. By assumption, $g\cou[o_i'] \leq o_i$ and $g\cou[a'] \in \clop X$ for all $a' \in \clop X'$. So, for any $a' \leq o_i'$, $g\cou[a'] \leq o_i$ and therefore $\clop_n g [I_{o_i'}'] \subseteq \{a \in \clop X \mid a \leq o_i\} = I_{o_i}$.

Then, as in the proof of Theorem \ref{dual}, we have two contravariant functors which are obviously faithful and whose two compositions are naturally isomorphic to the identity functors of the two categories. The assertion is proved.
\end{proof}

\begin{corollary}\label{mvn}
For all $n \in \omega$, the category $\MV_n$ of $n$-valued MV-algebras is dual to category $\ST_n$ of Stone spaces with $n-1$ distinguished open sets. 
\end{corollary}
\begin{proof}
By \cite[Theorem 22]{dnl} $\MV_n$ is equivalent to $\cat BR_n$ which, on its turn, is isomorphic to $\BB_n$ by Lemma \ref{relid}. Then the thesis follows from Theorem~\ref{stonen}.
\end{proof}

Now we shall combine our duality with the one presented in \cite{cidumu}. In that paper, the authors construct a category, denoted by $\C$, as follows.

First, consider the set $\N$ of positive natural numbers, equipped with the divisibility order $\mid$. $\la \N, \mid \ra$ is a distributive lattice, denoted by $\N_d$, with the join and meet of two numbers given by, respectively, the least common multiple and the greatest common divisor.
\begin{definition}\label{super}
A \emph{supernatural number} is a function $\nu: \P \lto \omega \cup \{\omega\}$, where $\P$ denotes the set of prime numbers. For any two supernatural numbers $\nu$ and $\mu$, we write $\nu \leq \mu$ iff $\nu(p) \leq \mu(p)$ for all $p \in \P$.
\end{definition}

Regarding $\nu$ as a list of exponents for the sequence of prime numbers, supernatural numbers can be seen as infinite formal products of (possibly infinite) powers of prime numbers. Then natural numbers can be identified with supernatural numbers with finite support and whose range is included in $\omega$, and the order relation of the supernaturals can be seen as the natural extension of the divisibility order in $\N$. With such an order, the supernatural numbers become a locale which will be denoted by $\G$.

The topology on $\G$ is defined as the one having as an open basis all sets of the form
$$U_n \bydef \{\nu \in \G : \nu > n \}, \textrm{ with } n \in \N.$$
By abuse of notation, $\G$ shall also denote the resulting topological space and it is worth remarking that the given topology coincides with the Scott topology \cite{scott}.

The objects of the category $\C$ of \emph{multisets} are pairs $\la \tau, s\ra$ such that $\tau = \la X, \Omega\ra$ is a Stone space and $s$ is a continuous functions from $\tau$ to $\G$. A morphism $f: \la \tau, s\ra \lto \la \tau', s'\ra$ is a continuous function from $\tau$ to $\tau'$ such that $\tau \geq \tau'\circ f$ with respect to the pointwise order.   

The main theorem of \cite{cidumu} can be stated as follows
\begin{theorem}\cite[Theorem 6.8]{cidumu}\label{cdmdual}
The category $\C$ is dual to the full subcategory of $\MV$ whose objects are locally finite MV-algebras.
\end{theorem}

Although not all locally finite MV-algebras are lcc, many lcc MV-algebras are locally finite (e.g., the liminary ones, as previously shown). So, taking into account also the results schematized in Table 1 of \cite{cidumu}, we can compose the duality of Theorem \ref{cdmdual} with some restrictions of our one, thus obtaining the situation described below.

Let $\approx$ mean ``categorically equivalent to'' and $\app$ mean ``dual to''. Moreover, for all $n < \omega$, let us denote by $\SMV(n)$ the subcategory of $\SMV$ whose objects' open sets have range included in \L$_n$ or, equivalently, whose single-point subspaces have clopen algebras isomorphic to some \L$_m$, with $m-1 \mid n-1$, and by $\SMV(\uparrow n)$ the union of all the $\SMV(k)$ for $k \leq n$. Then we have the following equivalences and dualities, where, from item 2 on, each category of MV-algebras is a subcategory of the one of the subsequent item.
\begin{enumerate}
\item Finite MV-algebras $\app$ Finite-valued finite Stone MV-spaces $\approx$ Finite multisets.
\item Boolean algebras $\app$ Classical Stone spaces $\approx$ Multisets with $s \equiv 1$.
\item $n$-Homogeneous MV-algebras (Boolean products of copies of \L$_n$) $\app$ $n$-valued Stone MV-spaces $\approx$ Multisets with $s \equiv n-1$.
\item $\MV_n \approx \cat BR_n \approx \BB_n \app \ST_n \approx \SMV(n) \approx$ Multisets with $s < n$. 
\item $n$-Bounded MV-algebras (Boolean products of \L$_k$, $2 \leq k \leq n$) $\app \SMV(\uparrow n) \approx$ Multisets with $s$ such that $\forall x \in X \ \exists k < n: s(x) = k$.
\item Liminary MV-algebras $\app$ Strongly compact $\SMV \approx$ Multisets with $s[X] \subseteq \omega$.
\end{enumerate}

As we already mentioned, MV-algebras boast a rather large number of topological dualities. Unfortunately, as far as we know, none of them (including the one presented here) covers the whole category of MV-algebras, and the scenario of the ``dualizable'' subcategories of $\MV$ is pretty wild. Nonetheless, we believe that our duality has essentially two distinguishing features. The first one is that it returns MV-algebras to their most natural --- yet too often neglected in the duality theory --- logical and set-theoretic environment: fuzziness. The second (and, probably, most important) feature is given by the fact that, although $\MVlcc$ is strictly contained in $\MVs$, Theorem \ref{lccfunc} gives the former a leading role among the countless subcategories of the latter, and represents, in our opinion, a further strong motivation for the study of limit cut complete MV-algebras.

\end{document}